\documentclass[11pt,english]{article}
\usepackage{graphicx} 
\usepackage{amssymb, amsthm, amsmath, fullpage, hyperref}
\usepackage{graphicx, enumerate}
\usepackage{tikz-network}
\usetikzlibrary{positioning, 
                quotes}
\usetikzlibrary{calc}
\usepackage[utf8]{inputenc}
\usepackage{cleveref}
\usepackage{amsfonts}
\usepackage{amssymb}
\usepackage{amsmath,amsthm,dsfont}
\usepackage[english]{babel}
\usepackage{makeidx}
\usepackage{mathtools}
\usepackage{bbm}
\usepackage{fullpage,amsfonts,amssymb,epsfig,epstopdf,amsmath,titling,url,array,titlesec}
\usepackage{mathtools}
\usepackage{authblk}
\usepackage{tikz}
\usetikzlibrary{trees, positioning, shapes, arrows.meta}
\tikzset{
  box/.style={draw, rectangle, rounded corners, minimum width=3.2cm, minimum height=0.9cm, align=center, font=\small},
  smallbox/.style={draw, rectangle, rounded corners, minimum width=2.5cm, minimum height=0.7cm, align=center, font=\small},
  line/.style={draw, -{Stealth[scale=1.1]}, thick},
}
\usepackage{float}
\usepackage{makeidx}
\usepackage[mathscr]{euscript}
\usepackage{nomencl}
\usepackage{comment}
\usepackage{titlesec}
\titleclass{\subsubsubsection}{straight}[\subsubsection]
\newcounter{subsubsubsection}[subsubsection]
\renewcommand\thesubsubsubsection{\thesubsubsection.\arabic{subsubsubsection}}
\titleformat{\subsubsubsection}
  {\normalfont\normalsize\bfseries}{\thesubsubsubsection}{1em}{}
  \titlespacing*{\subsubsubsection}
  {0pt}{3.25ex plus 1ex minus .2ex}{1.5ex plus .2ex}
\theoremstyle{plain}
\newtheorem{thm}{Theorem}[section]

\theoremstyle{definition}

\newtheorem{rem}[thm]{Remark}

\newtheorem{claim}[thm]{Claim}
\newtheorem*{notn}{Notation}
\newcommand\numberthis{\addtocounter{equation}{1}\tag{\theequation}}
\newtheorem{case}{Case}
\newtheorem{subcase}{Subcase}[case]
\newtheorem{subsubcase}{Subsubcase}[subcase]
\providecommand{\keywords}[1]
{
  \textbf{\textit{Keywords:}} #1
}
\allowdisplaybreaks

\author[1]{Rajat Adak}
\author[2]{L. Sunil Chandran}
\affil[1]{\texttt{rajatadak@iisc.ac.in}} 
\affil[2]{\texttt{sunil@iisc.ac.in}}
\affil[1,2]{Department of Computer Science and Automation, Indian Institute of Science, Bangalore, India}
\date{}
\title{Vertex-Based Localization of Tur\'{a}n's Theorem }
\begin{document}

\maketitle
\begin{abstract}
    Let $G$ be a simple graph with $n$ vertices and $m$ edges. According to Tur\'{a}n's theorem, if $G$ is $K_{r+1}$-free, then $m \leq |E(T(n, r))|,$ where $T(n, r)$ denotes the Tur\'{a}n graph on $n$ vertices with a maximum clique of order $r$. A limitation of this statement is that it does not give an expression in terms of $n$ and $r$. A widely used version of Tur\'{a}n's theorem states that for an $n$-vertex $K_{r+1}$-free graph,
$m \leq \left\lfloor \frac{n^2(r-1)}{2r} \right\rfloor.$
Though this bound is often more convenient, it is not the same as the original statement. In particular, the class of extremal graphs for this bound, say $\mathcal{S}$, is a proper subset of the set of Tur\'{a}n graphs. In this paper, we generalize this result as follows: For each $v \in V(G)$, let $c(v)$ be the order of the largest clique that contains $v$. We show that  
\[ m \leq \left\lfloor\frac{n}{2}\sum_{v\in V(G)}\frac{c(v)-1}{c(v)}\right\rfloor\]
Furthermore, we characterize the class of extremal graphs that attain equality in this bound. Interestingly, this class contains two extra non-Tur\'{a}n graphs other than the graphs in $\mathcal{S}$.
\end{abstract}
\section{Introduction}
Extremal graph theory studies the maximum or minimum values that certain graph parameters (typically the number of edges or vertices) can achieve subject to given constraints, and the corresponding extremal structures. Tur\'{a}n's Theorem is one of the most prominent results in this field.
\begin{notn}
    Let $T(n,r)$ denote the Tur\'{a}n Graph on $n$ vertices with a maximum clique of order $r$. It is a complete multipartite graph of order $n$ with $r$ classes of size either $\lceil\frac{n}{r}\rceil$ or $\lfloor\frac{n}{r}\rfloor$ in the partition. 
\end{notn}
\begin{thm}\label{th:Turan}\emph{(Tur\'{a}n \cite{turan1941egy})}
    For a simple $K_{r+1}$-free graph $G$ with $n$ vertices,
    \begin{equation*}
        |E(G)| \leq |E(T(n,r))|
    \end{equation*}
    and equality holds if and only if $G \cong T(n,r)$.
\end{thm}
\noindent Observe that the bound stated in \Cref{th:Turan} is not explicitly formulated as a function of 
$n$ and $r$. Consequently, a commonly used quantitative version of this bound is given as:
\begin{rem}\label{fturan_rem}
    For a simple $K_{r+1}$-free graph $G$ with $n$-vertices,
    \[|E(G)| \leq \left\lfloor \frac{n^2(r-1)}{2r} \right\rfloor\]
\end{rem}
\noindent The bound presented in \cref{fturan_rem} is not equivalent to the original formulation in \Cref{th:Turan}. In particular, the class of extremal graphs 
$\mathcal{S}$ corresponding to \cref{fturan_rem} forms a proper subset of the set of all Tur\'{a}n graphs.

\vspace{2mm}
\noindent\textbf{$\bullet$ Characterizing the set $\mathcal{S}$}
\newline Suppose that $G \in \mathcal{S}$ and $G \cong T(n,r)$ where $n =tr+s$, such that $0\leq s < r$. Thus, there are $s$ classes of size $(t+1)$ and $(r-s)$ classes of size $t$ in $G$. Therefore, we get:
\begin{equation}\label{no. of edges}
    |E(G)| = {n\choose2} - s{t+1\choose 2}-(r-s){t\choose 2} = \frac{r-1}{2r}(n^2-s^2) + {s\choose2}
\end{equation}
Since $G \in \mathcal{S}$, we also have;
\begin{equation}\label{G in S}
    |E(G)| = \left\lfloor \frac{n^2(r-1)}{2r} \right\rfloor \implies \frac{n^2(r-1)}{2r}  - |E(G)| <1
\end{equation}
Plugging in the value of $|E(G)|$ from \cref{no. of edges} into \cref{G in S} we get;
\begin{equation*}
    \begin{aligned}
        \frac{s^2(r-1)}{2r}-{s\choose2} <1 \implies \frac{s}{2r}(s(r-1)-r(s-1))< 1 \\ \implies \frac{s(r-s)}{2r}<1 \implies s(r-s) <2r \implies r(s-2) <s^2
    \end{aligned}
\end{equation*}
Note that, if $s \leq 2$ then $r(s-2) <s^2$ for any $r$. If $s>2$ we need $r < \frac{s^2}{s-2}$ for the condition to be true. Thus we can define;
\begin{equation}\label{defnS}
    \mathcal{S} = \{T(n,r) \mid \text{either } s \leq 2 \text{ or } r < \frac{s^2}{s-2}; \text{ where }n \equiv s\ (\text{mod } r)\}
\end{equation}
\begin{thm}
\label{fturan}
    For a simple $K_{r+1}$-free graph $G$ on $n$-vertices,
    \[|E(G)| \leq \left\lfloor \frac{n^2(r-1)}{2r} \right\rfloor\]  
    equality holds if and only if $G \cong T(n,r)$ and $G \in \mathcal{S}$.
\end{thm}
\noindent A \textit{relaxed} form of \Cref{fturan} is obtained by omitting the floor function from the right-hand side of the bound.

\begin{thm}\label{wturan}
For a simple $K_{r+1}$-free graph $G$ on $n$ vertices,
\[
|E(G)| \leq \frac{n^2(r-1)}{2r}.
\]
equality holds if and only if $G \cong T(n,r)$ and $r \mid n$; that is, $G$ is a regular Tur\'{a}n graph with $r$ equal-sized classes.
\end{thm}

\noindent
This simplification is often adopted by authors for convenience (even if it is a weaker version) because the characterization of the extremal graphs becomes very simple in the absence of the floor function. However, it is important to note that, since $|E(G)|$ must be an integer, the tightest possible bound in terms of $n$ and $r$ is provided by \Cref{fturan}. This is further underscored by the fact that the extremal graph class associated with \Cref{fturan} strictly contains that of \Cref{wturan}. 
 This is because, if $G$ is an extremal graph for \Cref{wturan}, then $n \equiv 0\ (\text{mod }r)$, that is, $s = 0$, thus $G \in \mathcal{S}$.
\subsection{Localization}
 Bradač \cite{bradac} gave a generalization of \Cref{wturan} by associating weights to the edges of the graph. Subsequently, Malec and Tompkins \cite{DBLP:journals/ejc/MalecT23} provided an alternative proof for this generalized
 version. In their work, the edge weights were defined as follows:
 \begin{equation*}
    k(e) = max\{r \mid e\ \text{occurs in a subgraph of $G$ isomorphic to}\ K_r\}
\end{equation*}
\begin{thm}\label{thm:GTuran}\emph{(Bradač \cite{bradac} and Malec-Tompkins \cite{DBLP:journals/ejc/MalecT23})} For a simple graph $G$ with $n$ vertices,
\begin{equation*}
    \sum_{e \in E(G)}\frac{k(e)}{k(e) -1} \leq \frac{n^2}{2}
\end{equation*}
and equality holds if and only if $G$ is regular Tur\'{a}n graph.
\end{thm}
\begin{rem}
    If $G$ is $K_{r+1}$ free, $k(e) \leq r$ for all $e \in E(G)$. Replacing $k(e)$ by $r$ in \Cref{thm:GTuran} yields the bound of \Cref{wturan} since $\frac{r}{r-1} \leq \frac{k(e)}{k(e)-1}$ for all $e \in E(G)$.
\end{rem}

\noindent Malec and Tompkins \cite{DBLP:journals/ejc/MalecT23} coined the term \textit{localization} to describe this type of generalization based on such weight assignments, as the weight of each edge depends solely on the local structures in which it is involved. Kirsch and Nir \cite{Kirsch_2024} extended \Cref{thm:GTuran} by assigning weights to cliques of any size, rather than just edges. Their result simultaneously extends a different extension of \Cref{wturan} due to Zykov \cite{Zykov1949-qr}. Arag\~{a}o and Souza \cite{aragao2024localised} further extended this by providing localization for the Graph Maclaurin Inequalities \cite{Khadzhiivanov1977-xg}.
\subsection{Vertex-Based Localization}
In Theorem~\ref{thm:GTuran}, localization is achieved by assigning weights to the edges of the graph. In contrast, we adopt the framework of \textit{vertex-based localization}, introduced in~\cite{adak}, wherein weights are assigned to the vertices rather than the edges. In this work, we present a vertex-based localization of \Cref{fturan}. Additionally, for the benefit of the interested reader, we include a concise proof of a vertex-based localization of \Cref{wturan} also. Although the bound in this weaker result follows as a consequence of the localized version of \Cref{fturan}, the extremal class in this case admits a simpler structural characterization. This allows for a direct comparison with \Cref{thm:GTuran}, concentrating only on the weaker formulation.

\section{Our Result}
Before presenting our main result, we introduce two specific graphs. Let $X \cong \overline{P_2 \cup P_3}$, commonly referred to in the literature as the \textit{Paraglider} graph, and let $Y \cong \overline{P_3}$, where $P_k$ denotes the path on $k$ vertices.

 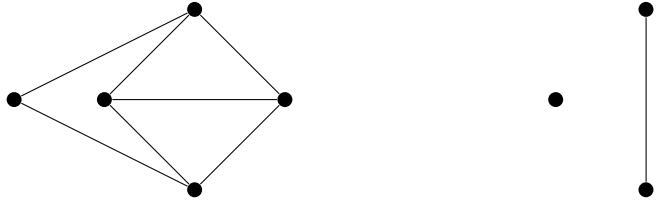
\begin{figure}[H]
    \centering
    \begin{tikzpicture}
    [scale=1.2, every node/.style={circle, fill=black, inner sep=2pt}]
    \centering
    \node (a) at (0,0) {};
    \node (b) at (2,1) {};
    \node (c) at (2,-1) {};
    \node (d) at (3,0) {};
    \node (e) at (1,0) {};

    \draw (a) -- (b);
    \draw (a) -- (c);

    \draw (b) -- (e);
    \draw (b) -- (d);
    \draw (c) -- (e);
    \draw (c) -- (d);
    \draw (e) -- (d);

    \node (f) at (7,1) {};
    \node (g) at (7,-1) {};
    \node (h) at (6,0) {};

    \draw (f) -- (g);

\end{tikzpicture}
\caption{Graph $X$ (left); Graph $Y$ (right)}
\end{figure}
\begin{notn}
    For a graph $G$, let $c(v)$ denote the weight of a vertex $v \in V(G)$, defined as;
\begin{equation*}
    c(v) = max\{r \mid v\ \text{occurs in a subgraph of $G$ isomorphic to}\ K_r\}
\end{equation*}
\end{notn}
\begin{thm}\label{th:main}
    For a simple graph $G$ with $n$ vertices
\begin{equation*}
    |E(G)| \leq \left\lfloor\frac{n}{2}\sum_{v\in V(G)}\frac{c(v)-1}{c(v)}\right\rfloor
\end{equation*}
Equality holds if and only if $G \in \mathcal{S} \cup \{X,Y\}$.
\end{thm} 
\subsection{Recovering \Cref{fturan}}
To recover the bound and the extremal graph for \Cref{fturan} from \Cref{th:main} we assume $G$ to be a $K_{r+1}$-free graph on $n$ vertices. Thus, $c(v) \leq r$, which implies $\frac{c(v)-1}{c(v)} \leq \frac{r-1}{r}$ for all $v \in VG)$. Therefore, we get:
\begin{equation}\label{eq2}
    |E(G)| \leq \left\lfloor\frac{n}{2}\sum_{v \in V(G)} \frac{c(v) -1}{c(v)}\right\rfloor \leq \left\lfloor\frac{n}{2}\sum_{v\in V(G)}\frac{r-1}{r}\right\rfloor = \left\lfloor\frac{n^2(r-1)}{2r}\right\rfloor
\end{equation}
Thus, we recover the bound stated in \Cref{fturan}. The first inequality in \cref{eq2} attains equality if and only if $G \in \mathcal{S} \cup \{X, Y\}$, from \Cref{th:main}. It is easy to verify that when $G \cong X$, we have $r = 3$, and when $G \cong Y$, we have $r = 2$; however, in both cases, the inequality is not tight. Thus, $G$ must be in $ \mathcal{S}$. That is, $G$ is a Tur\'{a}n graph, therefore $c(v) = r$ for all $v \in V(G)$. This gives equality in the second inequality of \cref{eq2}. Thus $G \cong T(n,r)$ and $G \in \mathcal{S}$.
\subsection{Proof of \Cref{th:main}}
\begin{proof}
The proof is based on induction over $|V(G)|$.
\vspace{2mm}
\newline\textit{Base Case:} If $|V(G)| = 1$, then $m = 0$, and the claim is trivially true since an isolated vertex has weight $1$, thus making the right-hand side $0$.
\vspace{1mm}
\newline\textit{Induction Hypothesis:} Suppose that the claim in \Cref{th:main} is true for any graph with less than $|V(G)|$ vertices.
\vspace{1mm}
\newline\textit{Induction Step:} Let $C$ be the largest clique in $G$. Let $u \in V(C)$. Clearly $c(u) = |V(C)|$. Suppose $G' = G \setminus V(C)$, that is $G'$ is obtained from $G$ by removing the vertices of $C$.
\subsubsection{Proof of Inequality}\label{Proof of Ineq}
We partition the edges of $G$ into three pairwise disjoint classes; $E(C), E(G'), \text{and } E(G': C)$, where $E(G':C)$ represents the set of edges having one endpoint in $G'$ and other in $C$. Thus we have;
\begin{equation}\label{eq3}
    |E(G)| = |E(C)| + |E(G')| + |E(G':C)|
\end{equation}
Now we count the number of edges in each class. Since $C$ is a maximum clique and of order $c(u)$, we have $c(v) = c(u)$ for all $v \in V(C)$. Therefore;
\begin{equation}\label{eq4}
    |E(C)| = \frac{c(u)(c(u) -1)}{2} =\frac{c(u)}{2}\sum_{v \in V(C)}\frac{c(v)-1}{c(v)} \in \mathbb{Z}
\end{equation}
Let $c_{G'}(v)$ denote the weight of a vertex $v \in V(G')$ restricted to $G'$. Thus $c_{G'}(v) \leq c(v)$ for all $v \in V(G')$. This implies, $\frac{c_{G'}(v) -1}{c_{G'}(v)} \leq \frac{c(v) -1}{c(v)}$ for all $v \in V(G')$. 
\newline $|V(G')| = n -|V(C)| = n - c(u)$, therefore by the induction hypothesis, we get the following;
\begin{equation}\label{ineq12}
    |E(G')| \leq \left\lfloor\frac{n-c(u)}{2}\sum_{v \in V(G')}\frac{c_{G'}(v) -1}{c_{G'}(v)}\right\rfloor \leq \left\lfloor\frac{n-c(u)}{2}\sum_{v \in V(G')}\frac{c(v) -1}{c(v)}\right\rfloor
\end{equation}
For all $v \in V(G')$, define $N_v = \{ w \in V(C) \mid w \in N(v)\}$. Clearly $N_v \cup \{v\}$ induces a clique for all $v \in V(G')$. Therefore, \begin{equation}\label{ineq3}
    |N_v| +1 \leq c(v)
\end{equation}
Note that the number of edges in $E(G':C)$ having $v \in V(G')$ as an endpoint is $|N_v|$. Thus using \cref{ineq3} we get the following;
\begin{equation}\label{eq6}
    |E(G':C)| = \sum_{v \in V(G')}|N_v| \leq \sum_{v\in V(G')}(c(v) -1) = \sum_{v\in V(G')}c(v)\frac{c(v) -1}{c(v)}\in \mathbb{Z} 
\end{equation}
From \cref{eq3,eq4,ineq12,eq6} we get;

  \begin{align}
    |E(G)| &\leq \frac{c(u)}{2}\sum_{v \in V(C)}\frac{c(v)-1}{c(v)} + \left\lfloor\frac{n-c(u)}{2}\sum_{v \in V(G')}\frac{c(v) -1}{c(v)}\right\rfloor +  \sum_{v\in V(G')}c(v)\frac{c(v) -1}{c(v)} \nonumber
    \\
    &= \left\lfloor\frac{c(u)}{2}\sum_{v \in V(C)}\frac{c(v)-1}{c(v)}\right\rfloor + \left\lfloor\frac{n-c(u)}{2}\sum_{v \in V(G')}\frac{c(v) -1}{c(v)}\right\rfloor +  \left\lfloor\sum_{v\in V(G')}c(v)\frac{c(v) -1}{c(v)} \right\rfloor\nonumber
    \\
        \text{\{Since,}& \text{ the first and the third terms are integers\}} \nonumber\\
    &\leq \left\lfloor\frac{c(u)}{2}\sum_{v \in V(C)}\frac{c(v)-1}{c(v)}\right\rfloor + \left\lfloor\frac{n-c(u)}{2}\sum_{v \in V(G')}\frac{c(v) -1}{c(v)}\right\rfloor +  \left\lfloor c(u)\sum_{v\in V(G')}\frac{c(v) -1}{c(v)} \right\rfloor\label{ineq4}
    \\
    &\leq \left\lfloor\frac{c(u)}{2}\sum_{v \in V(C)}\frac{c(v)-1}{c(v)}+ \frac{n-c(u)}{2}\sum_{v \in V(G')}\frac{c(v) -1}{c(v)} +   c(u)\sum_{v\in V(G')}\frac{c(v) -1}{c(v)} \right\rfloor\label{ineq5}
    \\
    &= \left\lfloor\frac{c(u)}{2}\sum_{v \in V(C)}\frac{c(v)-1}{c(v)}+ \frac{n}{2}\sum_{v \in V(G')}\frac{c(v) -1}{c(v)} +   \frac{c(u)}{2}\sum_{v\in V(G')}\frac{c(v) -1}{c(v)} \right\rfloor\nonumber
    \\
    & = \left\lfloor\frac{n}{2}\sum_{v \in V(G')}\frac{c(v) -1}{c(v)} +  \frac{c(u)}{2}\sum_{v \in V(G)}\frac{c(v)-1}{c(v)}\right\rfloor\nonumber
    \\
    & \leq \left\lfloor\frac{n}{2}\sum_{v \in V(G')}\frac{c(v) -1}{c(v)} +  \frac{c(u)}{2}\sum_{v \in V(G)}\frac{c(u)-1}{c(u)}\right\rfloor \label{ineq6}
    \\
    &=\left\lfloor\frac{n}{2}\sum_{v \in V(G')}\frac{c(v) -1}{c(v)} +  \frac{n}{2}(c(u)-1)\right\rfloor\nonumber \\
    &=\left\lfloor\frac{n}{2}\sum_{v \in V(G')}\frac{c(v) -1}{c(v)} + \frac{n}{2}\sum_{v \in V(C)}\frac{c(v) -1}{c(v)}\right\rfloor = \left\lfloor\frac{n}{2}\sum_{v\in V(G)}\frac{c(v)-1}{c(v)}\right\rfloor\nonumber
 \end{align}

Thus, we get the required bound as in \Cref{th:main}   
\end{proof} 
\subsubsection{Characterizing extremal graph class} 
\subsubsubsection{The only if part}\label{onlyif}
We will assume that equality holds in \Cref{th:main} and show that this implies $G \in \mathcal{S} \cup \{X,Y\}$. For equality in the bound of \Cref{th:main}, the inequalities in \cref{ineq12},(\ref{ineq3}),(\ref{ineq4}),(\ref{ineq5}), and (\ref{ineq6}) must become equalities. 

\noindent\newline For equality in the first inequality of \ref{ineq12}, using the induction hypothesis we get that $G' \in \mathcal{S}\cup \{X,Y\}$. For equality in the second inequality of \ref{ineq12}, we must have;
\begin{align}&\left\lfloor\frac{n-c(u)}{2}\sum_{v \in V(G')}\frac{c_{G'}(v) -1}{c_{G'}(v)}\right\rfloor = \left\lfloor\frac{n-c(u)}{2}\sum_{v \in V(G')}\frac{c(v) -1}{c(v)}\right\rfloor\label{1ineq2}\\
\implies&
 \left(\frac{n-c(u)}{2}\sum_{v \in V(G')}\frac{c(v) -1}{c(v)}\right)-\left(\frac{n-c(u)}{2}\sum_{v \in V(G')}\frac{c_{G'}(v) -1}{c_{G'}(v)}\right)<1\label{2ineq2}
\end{align}
\newline For equality in \ref{ineq3}, for all $v \in V(G')$ we must have;
\begin{equation}\label{1ineq3}
    |N_v| = c(v) -1
\end{equation} 
Note that, $\left\lfloor\sum_{v\in V(G')}c(v)\frac{c(v) -1}{c(v)} \right\rfloor$ is an integer. To get equality in \ref{ineq4} we need;
\begin{align}
 &\left\lfloor c(u)\sum_{v\in V(G')}\frac{c(v) -1}{c(v)} \right\rfloor =  \sum_{v\in V(G')}c(v)\frac{c(v) -1}{c(v)} \label{1ineq4}\\
 \implies& \left(c(u)\sum_{v\in V(G')}\frac{c(v) -1}{c(v)}  \right) - \left(\sum_{v\in V(G')}c(v)\frac{c(v) -1}{c(v)} \right) <1 \label{2ineq4}
\end{align}

\noindent\newline Let $frac\{x\} = x - \lfloor x\rfloor$. We get equality in \ref{ineq5} if and only if;
\begin{equation}
     frac\left\{\frac{n-c(u)}{2}\sum_{v \in V(G')}\frac{c(v) -1}{c(v)}\right\} +    frac\left\{c(u)\sum_{v\in V(G')}\frac{c(v) -1}{c(v)}\right\}<1\label{1ineq5}
\end{equation}
\newline Finally, for equality in \ref{ineq6} we must have;
\begin{align}
   &\left\lfloor\frac{n}{2}\sum_{v \in V(G')}\frac{c(v) -1}{c(v)} +  \frac{c(u)}{2}\sum_{v \in V(G)}\frac{c(v)-1}{c(v)}\right\rfloor =
 \left\lfloor\frac{n}{2}\sum_{v \in V(G')}\frac{c(v) -1}{c(v)} +  \frac{c(u)}{2}\sum_{v \in V(G)}\frac{c(u)-1}{c(u)}\right\rfloor\nonumber 
\end{align}

\vspace{2mm}

\noindent\textbf{$\bullet$ Preprocessing Steps}
\begin{rem}\label{notclique}
   If $G$ is a clique, then it is clearly extremal for \Cref{th:main}, and $G \in \mathcal{S}$. Therefore, we henceforth assume that $G$ is not a clique and thus $n - c(u) \neq 0$.

\end{rem}
\begin{claim}\label{G in S}
    If $G$ is an extremal graph for \Cref{th:main} and $G$ is a Tur\'{a}n Graph, then $G \in \mathcal{S}$.
    \begin{proof} To see this we first note that if $G$ is extremal for \Cref{th:main} then,
    \begin{equation*}
        |E(G)| = \left\lfloor\frac{n}{2}\sum_{v\in V(G)}\frac{c(v)-1}{c(v)}\right\rfloor
    \end{equation*}
    Now since $G$ is assumed to be a Tur\'{a}n graph, $c(v) = r$ for all $v\in V(G)$, where $r$ is the number of classes in $G$. Thus we get; 
    \begin{equation*}
        \left\lfloor\frac{n}{2}\sum_{v\in V(G)}\frac{c(v)-1}{c(v)}\right\rfloor = \left\lfloor\frac{n}{2}\sum_{v\in V(G)}\frac{r-1}{r}\right\rfloor = \left\lfloor\frac{n^2(r-1)}{2r}\right\rfloor
    \end{equation*}
    Combining the two observations above, we conclude that
\begin{equation*}
    |E(G)| = \left\lfloor\frac{n^2(r-1)}{2r}\right\rfloor \implies G \in \mathcal{S}.
\end{equation*}

    \end{proof}
\end{claim}
\begin{claim}\label{conn}
    If $G$ is an extremal graph for \Cref{th:main}, and $G$ is neither a collection of isolated vertices (edgeless graph) nor $Y$, then $G$ is a connected graph.
    \begin{proof}
        Suppose $G$ is a disconnected extremal graph and $G$ is neither an edgeless graph nor the graph
        $Y$. If $G$ contains two non-trivial components, we can add an edge connecting these components without increasing the weight of any vertex but increasing the edge count. Thus, $G$ contains exactly one non-trivial component, say $H$. Clearly, $c(v) \geq 2$ for all $v \in V(H)$. Let $x \in V(G)\setminus V(H)$. Clearly, $c(x) =1$. Let $G' = G \setminus\{x\}$. Note that $|E(G)| = |E(G')|$. Since $G$ is extremal;
        \begin{equation}\label{connectivity1}
            |E(G')| = |E(G)| = \left\lfloor\frac{n}{2}\sum_{v \in V(G)}\frac{c(v)-1}{c(v)}\right\rfloor
        \end{equation}
        After the removal of vertex $x$ from $G$, the weights of the remaining vertices remain the same. Now applying the bound of \Cref{th:main} on $G'$, then using $c(x) =1$, and from \cref{connectivity1} we get;
        \begin{equation}\label{connectivity2}
            |E(G')| \leq \left\lfloor\frac{n-1}{2}\sum_{v \in V(G')}\frac{c(v)-1}{c(v)}\right\rfloor = \left\lfloor\frac{n-1}{2}\sum_{v \in V(G)}\frac{c(v)-1}{c(v)}\right\rfloor \leq \left\lfloor\frac{n}{2}\sum_{v \in V(G)}\frac{c(v)-1}{c(v)}\right\rfloor = |E(G')|
        \end{equation}
        Therefore from \cref{connectivity2} we must have;
        
        \begin{align}
            &\left\lfloor\frac{n}{2}\sum_{v \in V(G)}\frac{c(v)-1}{c(v)}\right\rfloor = \left\lfloor\frac{n-1}{2}\sum_{v \in V(G)}\frac{c(v)-1}{c(v)}\right\rfloor\nonumber\\
            \implies& \left(\frac{n}{2}\sum_{v \in V(G)}\frac{c(v)-1}{c(v)}\right)-\left(\frac{n-1}{2}\sum_{v \in V(G)}\frac{c(v)-1}{c(v)}\right)<1\nonumber\\
           \implies& \frac{1}{2}\sum_{v \in V(G)}\frac{c(v)-1}{c(v)}<1 
           \implies \sum_{v \in V(G)}\frac{c(v)-1}{c(v)}<2\label{connectivity3}
        \end{align}
    \end{proof}
    \noindent Since $\frac{c(v)-1}{c(v)}\geq \frac{1}{2}$ for all $v \in V(H)$. Thus using \cref{connectivity3}, we get $|V(H)| <4$. Therefore, $|V(H)| \in \{2,3\}$. Let the number of isolated vertices in $G$ be $k\geq 1.$
    \newline Suppose $|V(H)| = 2$, then $|E(G)| = 1$. Then we get;
    \[1= \left\lfloor\frac{k+2}{2}\left(\frac{1}{2}+\frac{1}{2}\right)\right\rfloor =\left\lfloor\frac{k+2}{2}\right\rfloor \implies k \leq 1 \implies k =1\implies G\cong Y\]
    But we assumed that $G$ is not isomorphic to $Y$, therefore $|V(G)| = 3$. Then $H$ can either be a $K_3$ or a $P_3$. If $H \cong K_3$, then $|E(G)| = 3$ and $c(v) = 3$ for all $v \in V(H)$. Then we get;
    \[3 = \left\lfloor\frac{k+3}{2}\left(3\times \frac{2}{3}\right)\right\rfloor =\left\lfloor k+3\right\rfloor \implies k <1\]
    Now assuming $H \cong P_3$, $|E(G)| = 2$ and $c(v) =2$ for all $v \in V(H)$. Therefore we get;
    \[2 = \left\lfloor\frac{k+3}{2}\left(3\times \frac{1}{2}\right)\right\rfloor =\left\lfloor\frac{3(k+3)}{4}\right\rfloor \implies k <1\]
\noindent In both cases, we get a contradiction. Therefore, $G$ can not be disconnected.
    
\end{claim}
\begin{claim}\label{light}
    Let $G$ be an extremal graph for \Cref{th:main}. Let $C$ be a clique with the largest order in $G$, and $G' = G \setminus V(C)$, and $u \in V(C)$. Then one of the following holds;
    \begin{itemize}
        \item $c(v) = c(u)$ for all $v \in V(G)$.
        \item There exists $v_0 \in V(G)$ such that $c(v_0) = (c(u) -1)$ and $c(v) = c(u)$ for all $v \in V(G) \setminus\{v_0\}$.
    \end{itemize}
    \begin{proof}
        If $G$ is either an edgeless graph or the graph $Y$, then the claim is clearly true. Thus, using \cref{conn} we can assume $G$ to be connected and thus $c(v) \geq 2$ for all $v \in V(G)$. Let $L = \{ v \in V(G) \mid c(v) \neq c(u)\}$. Suppose that $L$ is non-empty. We know that $c(v) \leq c(u)$ for all $v \in V(G)$. Therefore, $c(v) < c(u)$ for all $v \in L$. Clearly, $L \subseteq V(G')$. Now from \cref{2ineq4} we get;
        \begin{align}
            1 &> \left(c(u)\sum_{v\in V(G')}\frac{c(v) -1}{c(v)}  \right) - \left(\sum_{v\in V(G')}c(v)\frac{c(v) -1}{c(v)} \right) \nonumber\\
            &= \left(c(u)\sum_{v\in V(G')}\frac{c(v) -1}{c(v)}  \right) - \left(\sum_{v\in V(G')\setminus L}c(v)\frac{c(v) -1}{c(v)}\right) - \left(\sum_{v\in L}c(v)\frac{c(v) -1}{c(v)}\right)\nonumber\\
            &\geq \left(c(u)\sum_{v\in V(G')}\frac{c(v) -1}{c(v)}  \right) - \left(c(u)\sum_{v\in V(G')\setminus L}\frac{c(v) -1}{c(v)}\right) - \left(\sum_{v\in L}c(v)\frac{c(v) -1}{c(v)}\right)\nonumber\\
            &=\left(\sum_{v\in L}(c(u) -c(v))\frac{c(v)-1}{c(v)}\right)\label{light1}
        \end{align}
        Since $c(v) \geq 2$, we get $\frac{c(v)-1}{c(v)} \geq \frac{1}{2}$ for all $v \in L$. Since for all $v \in L$, $c(v) < c(u)$, we get $c(u)-c(v) \geq 1$. Therefore, from \cref{light1} we get;
        \[1 > \left(\sum_{v\in L}(c(u) -c(v))\frac{c(v)-1}{c(v)}\right) \geq \sum_{v \in L}\frac{c(v)-1}{c(v)} \geq \frac{|L|}{2} \implies 2 > |L|\]
        Thus $|L|= 1$ and let $L = \{v_0\}$. Again from \cref{light1} we get; \[1 > (c(u)-c(v_0))\frac{c(v_0)-1}{c(v_0)}\geq \frac{c(u) - c(v_0)}{2} \implies 2 > c(u) - c(v_0)\]
        Thus $c(v_0) = c(u) -1$ and $c(v) = c(u)$ for all $v \in V(G) \setminus \{v_0\}$
        \end{proof}
\end{claim}

\noindent We now proceed with a case analysis based on whether $G' \in \mathcal{S} $, $ G' \cong X $, or $ G' \cong Y $. In each case, as guided by \cref{light}, we further consider two subcases depending on whether all vertices have equal weights or not. To aid the reader in following the structure of our argument, we provide a case tree for the more sophisticated Case 1, $G' \in \mathcal{S}$, which illustrates the progression of cases and subcases.

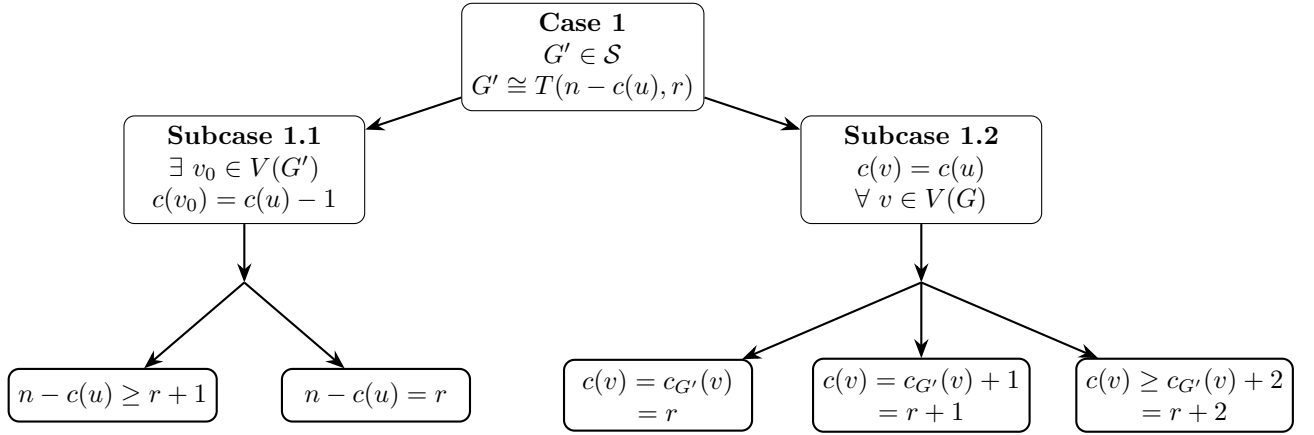
\begin{center}
\begin{tikzpicture}[
  level distance=1.5cm,
  edge from parent/.style={line},
  every node/.style={box},
  grow'=down,
  level 1/.style={sibling distance=90mm}, 
]

\node {\textbf{Case 1}\\$G' \in \mathcal{S}$\\$G' \cong T(n-c(u) ,r)$}
  child {node {\textbf{Subcase 1.2}\\$c(v) = c(u)$\\$\forall\ v \in V(G)$}
    child[grow'=down, level distance=1.5cm, sibling distance=35mm] {
      node[coordinate] {}
      child {node[smallbox] {$c(v) \geq c_{G'}(v) + 2$ \\ $= r+2$}}
      child {node[smallbox] {$c(v) = c_{G'}(v) + 1$ \\ $=r+1$}}
      child {node[smallbox] {$c(v) = c_{G'}(v)$ \\ $=r$}}
    }
  }
    child {node {\textbf{Subcase 1.1}\\$\exists$ $v_0 \in V(G')$\\$c(v_0) = c(u) - 1$}
    child[grow'=down, level distance=1.5cm, sibling distance=35mm] {
      node[coordinate] {}
      child {node[smallbox] {$n - c(u) = r$}}
      child {node[smallbox] {$n - c(u) \geq r+1$}}
    }
  };

\end{tikzpicture}
\end{center}
\vspace{2mm}

\begin{case} $G' \in \mathcal{S}$. Let $G' = T(n-c(u),r)$ for some $r\geq 1$ (since $G$ is not a clique, from \cref{notclique}). Therefore, $c_{G'}(v) = r$ for all $v \in V(G')$.
\begin{subcase}\label{sbc1.1} There exists $v_0 \in V(G')$ such that $c(v_0) = c(u) -1$ and $c(v) = c(u)$ for all $v \in V(G) \setminus\{v_0\}$.
\vspace{2mm}
\newline From \cref{1ineq5} we have;
\begin{equation}\label{1.11}
     frac\left\{\frac{n-c(u)}{2}\sum_{v \in V(G')}\frac{c(v) -1}{c(v)}\right\} +    frac\left\{c(u)\sum_{v\in V(G')}\frac{c(v) -1}{c(v)}\right\}<1
\end{equation}
Since $c(v) \geq c_{G'}(v)$ for all $v \in V(G')$, we get that $c(v_0) \geq c_{G'}(v_0) = r$. Therefore, for all $v \in V(G) \setminus\{v_0\}$, $c(v) = c(v_0) +1 \geq r+1$. Note that $(*)$ below, follows from \cref{1ineq2}. Thus we get;
\begin{align*}
    frac\left\{\frac{n-c(u)}{2} \sum_{v \in V(G')}\frac{c(v) -1}{c(v)}\right\} &= \frac{n-c(u)}{2}\sum_{v \in V(G')}\frac{c(v) -1}{c(v)} - \left\lfloor\frac{n-c(u)}{2}\sum_{v \in V(G')}\frac{c(v) -1}{c(v)}\right\rfloor\\
    &=\frac{n-c(u)}{2}\sum_{v \in V(G')}\frac{c(v) -1}{c(v)} - \left\lfloor\frac{n-c(u)}{2}\sum_{v \in V(G')}\frac{c_{G'}(v) -1}{c_{G'}(v)}\right\rfloor\ \ (*)\\
    &\geq \frac{n-c(u)}{2}\sum_{v \in V(G')}\frac{c(v) -1}{c(v)} - \left(\frac{n-c(u)}{2}\sum_{v \in V(G')}\frac{c_{G'}(v) -1}{c_{G'}(v)}\right)\\
    \geq \frac{n-c(u)}{2}&\left(\sum_{v\in V(G') \setminus\{v_0\}} \frac{r}{r+1}\right) + \frac{n-c(u)}{2}\left(\frac{r-1}{r}\right) - \left(\frac{n-c(u)}{2}\sum_{v\in V(G')} \frac{r-1}{r}\right)\\
    &=\frac{n-c(u)}{2}\sum_{v\in V(G') \setminus\{v_0\}} \left(\frac{r}{r+1} - \frac{r-1}{r}\right)\\
    &=\frac{(n-c(u)-1)(n-c(u))}{2r(r+1)}\numberthis\label{1.12}
\end{align*}
Note that $(**)$ follows from \cref{1ineq4}. Therefore we get;
\begin{align*}
    frac\left\{c(u)\sum_{v\in V(G')}\frac{c(v) -1}{c(v)}\right\}&=c(u)\sum_{v\in V(G')}\frac{c(v) -1}{c(v)} - \left\lfloor c(u)\sum_{v\in V(G')}\frac{c(v) -1}{c(v)}\right\rfloor\\
    &=c(u)\sum_{v\in V(G')}\frac{c(v) -1}{c(v)} - \left( \sum_{v\in V(G')}c(v)\frac{c(v) -1}{c(v)}\right)\ \ (**)\\
    &=(c(u)-c(v_0))\frac{c(v_0)-1}{c(v_0)} = \frac{c(v_0)-1}{c(v_0)}\numberthis\label{1.13}
\end{align*}
The last step follows from the assumption of \cref{sbc1.1}. Note that $G'$ is $r$-partite, therefore $n-c(u) \geq r$. Now substituting the bound from \cref{1.11,1.12} in \cref{1.13} we get;
\begin{align*}
    &\frac{(n-c(u)-1)(n-c(u))}{2r(r+1)} + \frac{c(v_0)-1}{c(v_0)} <1\\
    \implies& \frac{(n-c(u)-1)(n-c(u))}{2r(r+1)} < \frac{1}{c(v_0)} \leq \frac{1}{r}\\
    \implies& \frac{(r-1)r}{2r(r+1)} <\frac{1}{r} \implies r^2 - 3r -2 <0\\
    \implies& r <\frac{3 +\sqrt{17}}{2} \implies r \leq 3
\end{align*}
Thus, the possible values of $r$ are $1,2$ and $3$.

\begin{subsubcase} $n-c(u) \geq r+1$.
\newline From the previous analysis, we know that;
\begin{align*}
    &\frac{(n-c(u)-1)(n-c(u))}{2r(r+1)} < \frac{1}{r}\\
    \implies& \frac{r(r+1)}{2r(r+1)}<\frac{1}{r}\implies r <2
\end{align*}
Therefore, the only possibility for $r$ is $1$. Thus $c_{G'}(v) = 1$ for all $v \in V(G')$ and $|V(G')| = n -c(u) \geq 2$. Note that, $|V(C)| = c(u) \geq r+1 =2$. Since $G$ contains a clique of order at least $2$, using \cref{conn} we can say that $G$ is connected unless $G \cong Y$. Since  $n-c(u)\geq 2$, $G \not\cong Y$. Thus $G$ must be connected, which implies $c(v) \geq 2$ for all $v \in V(G)$. Therefore we get;
\[\left\lfloor\frac{n-c(u)}{2}\sum_{v \in V(G')}\frac{c(v)-1}{c(v)}\right\rfloor \geq \left\lfloor\frac{2}{2}\left(\frac{1}{2}+\frac{1}{2}\right)\right\rfloor = 1 \]
But on the other hand we have;
\[\left\lfloor\frac{n-c(u)}{2}\sum_{v \in V(G')}\frac{c_{G'}(v)-1}{c_{G'}(v)}\right\rfloor = 0\]
Therefore, \cref{1ineq2} is not satisfied, a contradiction.
\end{subsubcase}

\begin{subsubcase}$n-c(u) = r$
\begin{itemize}
    \item If $r =1 \implies c_{G'}(v) =1 \implies |V(C)| =c(u) \geq r+1 = 2$. Note that, in this case $V(G') = \{v_0\}$. If $|V(C)| =c(u) =2$, then $c(v_0) = c(u) -1 =1$ and therefore $v_0$ is isolated in $G$ itself. Thus, by \cref{conn}, clearly $G \cong Y$.
    \begin{figure}[H]
    \centering
    \begin{tikzpicture}
    [scale=1.2, every node/.style={circle, fill=black, inner sep=2pt}]
    \centering
    \node (f) at (7,1) {};
    \node (g) at (7,-1) {};
    \node (h) at (6,0) [draw=black, fill=none, circle, inner sep=2pt] {{$v_0$}};

    \draw (f) -- (g);

\end{tikzpicture}
\end{figure}
    Otherwise, $|V(C)| >2$. Note that $V(G') = \{v_0\}$ and $c(v_0) = c(u) -1$. From \cref{1ineq3} we get $|N_{v_0}| = c(v_0)-1 = c(u) -2$. Therefore, \[|E(G)| = |E(C)| + |N_{v_0}| = \frac{c(u)(c(u)-1)}{2} + (c(u)-2)\]
    \begin{align*}
        \frac{n}{2}\sum_{v \in V(G)}\frac{c(v)-1}{c(v)} &= \frac{(c(u)+1)}{2}\left(c(u) \frac{c(u) -1}{c(u)} + \frac{c(u)-2}{c(u)-1}\right)\\
        &=\frac{c(u)}{2}(c(u)-1) + \frac{c(u)-1}{2}+ \frac{(c(u)-1)(c(u)-2)}{2(c(u)-1)} + \frac{2(c(u)-2)}{2(c(u) -1)}
    \end{align*}
    We know that for $G$ to be extremal we need,
    \begin{align*}
        & \left\lfloor\frac{n}{2}\sum_{v \in V(G)}\frac{c(v)-1}{c(v)}\right\rfloor = |E(G)|\\
        \implies&\frac{n}{2}\sum_{v \in V(G)}\frac{c(v)-1}{c(v)} - |E(G)| <1\\
        \implies& \frac{c(u)-1}{2}-\frac{c(u)-2}{2} +\frac{c(u)-2}{c(u)-1}<1\\
        \implies& \frac{c(u)-2}{c(u)-1} < \frac{1}{2} \implies c(u) <3
    \end{align*}
    This contradicts the fact that $|V(C)| = c(u) >2$.
    \item If $r =2 \implies c_{G'}(v) =2$. Thus $G' = T(2,2)$, is a single edge. Let $t = c(u) \geq r+1 = 3$. Let $V(G') = \{v_0,v_1\}$ and $c(v_0) = t-1$.
    \newline Substituting $c(v_0) = t-1$ and $c(v_1) = t$ in \cref{1ineq5} we get that;
    \begin{align*}
        &frac\left\{\frac{t-1}{t} +\frac{t-2}{t-1}\right\}
        + frac\left\{\frac{t(t-1)}{t} + \frac{t(t-2)}{t-1}\right\}<1\\
        \implies& \left(1 - \frac{1}{t} - \frac{1}{t-1}\right) + \left(1 - \frac{1}{t-1}\right) <1\\
        \implies& 1 < \frac{1}{t} + \frac{2}{t-1} \implies t <4
    \end{align*}
    But $t \geq 3$, therefore $t = 3$. $|N_{v_0}| = c(v_0) - 1 = t-2 =1$ and $|N(v_1)| = c(v_1) - 1 = t -1 =2$, a contradiction. Note that $N_{v_0} \cap N_{v_1} = \emptyset $, otherwise we will get a $K_3$ containing $v_0, v_1$ and the common neighbor. Thus $c(v_0)$ will become $3>t-1$, a contradiction. Now $C$ is a $K_3$; $v_0$ is adjacent to one vertex of $C$ and $v_1$ is adjacent to the other two vertices of $C$. Therefore, $G \cong X$, the paraglider.
    
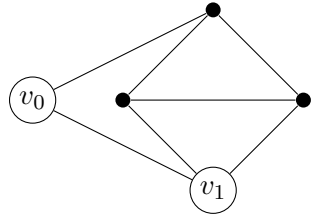
\begin{figure}[H]
    \centering
    \begin{tikzpicture}
    [scale=1.2, every node/.style={circle, fill=black, inner sep=2pt}]
    \centering
    \node (a) at (0,0) [draw=black, fill=none, circle, inner sep=2pt] {{$v_0$}};
    \node (b) at (2,1) {};
    \node (c) at (2,-1) [draw=black, fill=none, circle, inner sep=2pt] {{$v_1$}};
    \node (d) at (3,0) {};
    \node (e) at (1,0) {};

    \draw (a) -- (b);
    \draw (a) -- (c);

    \draw (b) -- (e);
    \draw (b) -- (d);
    \draw (c) -- (e);
    \draw (c) -- (d);
    \draw (e) -- (d);
\end{tikzpicture}
\end{figure}
\item If $r =3 \implies c_{G'}(v) =3$ Therefore $G' = T(3,3) = K_3$. Let $t = c(u) \geq r+1 = 4$. $v_0 \in V(G')$ and $c(v_0) = t-1$. From \cref{1ineq3}, $|N_{v_0}| = t-2$ and $|N_{v}| = t-1$ for all $v \in V(G) \setminus \{v_0\}$. Therefore;
\begin{equation}\label{*}
    |E(G)| = \frac{t(t-1)}{2} + 3 + 2(t-1) + (t-2) = \frac{t^4 + 4t^3 - 7t^2+2t}{2t(t-1)}
\end{equation}
\begin{equation}\label{**}
    \frac{n}{2}\sum_{v\in V(G)}\frac{c(v) -1}{c(v)} = \frac{t+3}{2}\left((t+2)\frac{t-1}{t}+\frac{t-2}{t-1}\right) = \frac{t^4+4t^3-2t^2-13t+6}{2t(t-1)}
\end{equation} 
If $G$ is extremal, then we have;
\begin{align*}
    &\frac{n}{2}\sum_{v\in V(G)}\frac{c(v) -1}{c(v)} - |E(G)| <1
\end{align*} 
Now substituting the values obtained in \cref{*} and \cref{**} we get;
\begin{align*}
    \implies& \frac{5t^2-15t+6}{2t(t-1)}<1\\
    \implies& 5t^2-15t+6 < 2t^2 -2t\\
    \implies& 3t^2-13t + 6<0 \implies t < \frac{13 + \sqrt{97}}{6} \implies t <4
\end{align*}
\newline This contradicts the fact that $t \geq 4$. Thus, there is no extremal graph, in this case.
\end{itemize}
Thus for $r= 1$ we got $G \cong Y$, for $r = 2$, $G \cong X$ and no extremal graph for $r =3$.
\vspace{2mm}
\end{subsubcase}
\end{subcase}
\begin{subcase}
$c(v) = c(u)$ for all $v \in V(G)$.
\begin{subsubcase}\label{1.2.1} $c(v) \geq c_{G'}(v)+2 \implies c(v) \geq r+2$.
\vspace{2mm}
\newline From \cref{2ineq2} we get;
\begin{align*}
    1&>\left(\frac{n-c(u)}{2}\sum_{v \in V(G')}\frac{c(v) -1}{c(v)}\right)-\left(\frac{n-c(u)}{2}\sum_{v \in V(G')}\frac{c_{G'}(v) -1}{c_{G'}(v)}\right)\\
    &\geq\left(\frac{n-c(u)}{2}\sum_{v \in V(G')}\frac{r +1}{r+2}\right)-\left(\frac{n-c(u)}{2}\sum_{v \in V(G')}\frac{r -1}{r}\right)\\
    &=\frac{n-c(u)}{2}\sum_{v\in V(G')}\left(\frac{r+1}{r+2} -\frac{r-1}{r}\right)\\
    &=\frac{(n-c(u))^2}{2}\frac{2}{r(r+2)} \implies (n-c(u))^2 <r(r+2)<(r+1)^2\\
    &\implies (n-c(u)) < r+1 \implies (n-c(u)) = r
\end{align*}
Thus, $G'= T(r,r)$ is a clique of size $r$. From \cref{1ineq3}we have $|N_v| = c(v) -1 = c(u) -1$, for all $v \in V(G')$. Therefore, all the vertices of $G'$ have $c(u)-1$ neighbors on $C$, a clique of size $c(u)$. Note that for $v_1,v_2 \in V(G')$, such that $v_1 \neq v_2$, $N_{v_1} \neq N_{v_2}$; otherwise $N_{v_1}\cup\{v_1,v_2\}$ will induce a $(c(u) +1)$-clique. Since each $v_i$ has exactly $|C|-1$ neighbors in $C$, and $N_{v_i} \neq N_{v_j}$ for $i\neq j$, each vertex $v_i \in V(G')$ has a unique non-neighbor in $C$. Without loss of generality let $u_i$ be the non-neighbor of $v_i$, for all $ i \in [r]$. Therefore, $G$ is a Tur\'{a}n graph, with classes, $\{v_1,u_1\},\dots,\{v_r,u_r\},\{u_{r+1}\},\dots,\{u_{c(u)}\}$. That is $G \cong T(c(u) + r, c(u))$. Thus, from \cref{G in S}, $G \in \mathcal{S}$.
\end{subsubcase}
\begin{subsubcase} $c(v) = c_{G'}(v) +1 \implies c(v) = r+1$.
\vspace{2mm}
\newline From \cref{2ineq2} we get that;
\begin{align*}
    1&>\left(\frac{n-c(u)}{2}\sum_{v \in V(G')}\frac{c(v) -1}{c(v)}\right)-\left(\frac{n-c(u)}{2}\sum_{v \in V(G')}\frac{c_{G'}(v) -1}{c_{G'}(v)}\right)\\
    &=\left(\frac{n-c(u)}{2}\sum_{v \in V(G')}\frac{r}{r+1}\right)-\left(\frac{n-c(u)}{2}\sum_{v \in V(G')}\frac{r -1}{r}\right)\\
    &=\frac{n-c(u)}{2}\sum_{v\in V(G')}\left(\frac{r}{r+1} -\frac{r-1}{r}\right)\\
    &=\frac{(n-c(u))^2}{2}\frac{1}{r(r+1)} \implies (n-c(u))^2 <2r(r+1)\\
\end{align*}
Suppose $n - c(u) = r +k$, where $k \geq 0$. Since $(r+k)^2 < 2r(r+1) \implies k \leq r$
\newline Since $G'$ is Tur\'{a}n, $G'$ comprises of $k$ many doubleton and $r-k$ many singleton parts. Therefore,
\[|E(G')| = \frac{1}{2}(2k(r+k-2)+(r-k)(r+k-1))\]
Since, for all $v \in V(G')$, $|N_v| = c(v) -1 = r$, $|E(G':C)| = r(r+k)$. Therefore;
\begin{align}\label{@}
    |E(G)|&= \frac{1}{2}(2k(r+k-2)+(r-k)(r+k-1)) + r(r+k) + \frac{r(r+1)}{2}\nonumber\\ &= 
\frac{4r^3 + 4kr^2 + 4r^2 + kr + k^2r + k^2 - 3k}{2(r+1)}
\end{align}
\begin{align}\label{@@}
    \frac{n}{2}\sum_{v \in V(G)}\frac{c(v)-1}{c(v)} = \frac{(2r+k+1)^2}{2}\left(\frac{r}{r+1}\right)  = \frac{4r^3 + 4r^2k + 4r^2 + rk^2 + 2rk + r}{2(r + 1)}
\end{align}
Since $G$ is extremal, substituting the values from \cref{@} and \cref{@@} we get;
\begin{align*}
    1>&\frac{n}{2}\sum_{v \in V(G)}\frac{c(v)-1}{c(v)} - |E(G)|=\frac{rk + r - k^2 + 3k}{2(r + 1)}\\
    \implies&2(r+1) >  rk + r - k^2 + 3k \implies k^2 - k(r+3) +(r+2) >0\\
    \implies&(k-(r+2))(k-1)>0
\end{align*}
If $k \geq 1 \implies k >(r+2)$. This contradicts the fact that $k \leq r$. Therefore $k = 0$. Thus $n-c(u) =r$ and $G' = T(r,r) = K_r$. Thus, following a similar argument as in \cref{1.2.1}, we get that each vertex $v \in V(G')$ has $|N_v| = |C| -1$ and an exclusive non-adjacent vertex in $C$. From this we infer as in \cref{1.2.1} that, $G\cong T(c(u) + r, c(u)) = T(2r+1, r+1)$. Since $G$ is Tur\'{a}n, from \cref{G in S} we get $G \in \mathcal{S}$.
\end{subsubcase}
\begin{subsubcase}\label{1.2.3}$c(v) = c_{G'}(v) \implies c(v) = r$.
\vspace{2mm}
\newline Recall that $G' = T(n-c(u),r)$. Let the classes of $G'$ be $P_1, P_2,\dots,P_r$, we have $||P_i| - |P_j||\leq 1$, where $i,j \in[r]$. Note that from \cref{1ineq3}, $|N_v| = c(v)-1 = r-1$, for all $v \in V(G')$. Let $v_i \in P_i$ and $v_j \in P_j$ where $i \neq j$, then $v_i$ is adjacent to $v_j$. Clearly as in the previous case, $N_{v_i} \neq N_{v_j}$, otherwise we get a clique of size $(c(u)+1)$.

\noindent\newline  Since for $v \in V(G')$, $|N_v| = r-1$ and $N_v \subseteq V(C)$, there are only ${r\choose r-1} = r$ possible choices for $N_v$. Pick a set of $r$ vertices say, $\{v_1, v_2, \dots, v_r\}$ such that $v _i \in P_i$ for all $i \in [r]$. We know that $N_{v_i} \neq N_{v_2}$, when $ i\neq j$. Since each $v_i$ must have an exclusive $N_{v_i}$, all the $r$ choices of $N_v$ are exhausted. It follows that for $x, y \in P_i$, $N_x = N_y$; otherwise we will get some $j \neq i$ and $z \in P_j$ such that either $N_x = N_z$ or $N_y = N_z$, a contradiction.

\noindent\newline Let $V(C) \setminus N_{v_i} = \{u_i\}$, and define $P'_i = P_i \cup \{u_i\}$. Note that $u_i \notin N(v)$ for all $v \in P_i$, but $u_i \in N(v)$ for all $v \in P_j$, and $u_i$ is adjacent to $u_j$ (since $u_i,u_j \in V(C)$), where $j \neq i$. Also $||P'_i| - |P'_j|| = ||P_i|+1-|P_j|-1| =||P_i|-|P_j|| \leq 1$. Therefore, $G$ is a Tur\'{a}n graph with classes $P'_1,P'_2,\dots,P'_r$. Thus from \cref{G in S} we get $G \in \mathcal{S}$.
\end{subsubcase}
\end{subcase}
\end{case}
\begin{case}
    $G' \cong X$
    \newline Since $X$ contains a $K_3$, clearly $c(u) \geq 3$. Note that for $4$ vertices of of $G \cong X$, $c_{G'}(v) = 3$ and for the remaining one vertex $c_{G'}(v) =2$.
\begin{subcase}
There exists $v_0 \in V(G')$ such that $c(v_0) = c(u) -1$ and $c(v) = c(u)$ for all $v \in V(G) \setminus\{v_0\}$.    
\vspace{2mm} 
\newline Suppose $c(u) \geq 4$, then $c(v_0) \geq 3$. From \cref{2ineq2} we get;
\begin{align*}
    1&> \frac{5}{2}\sum_{v \in V(G')}\frac{c(v)-1}{c(v)} - \frac{5}{2}\sum_{v \in V(G')}\frac{c_{G'}(v)-1}{c_{G'}(v)}\\
    &\geq \frac{5}{2}\left(4 \times \frac{3}{4} + \frac{2}{3}\right)- \frac{5}{2}\left(4\times\frac{2}{3}+\frac{1}{2}\right) = \frac{5}{4}
\end{align*}
We get a contradiction. Therefore $c(u)$ must be $3$. Thus $|N_v| = 2$ for all $v \in V(G')\setminus \{v_0\}$ and $|N_{v_0}| = 1$. We get;
\[|E(G)| = |E(G')| + |E(C)| + |E(G':C)| = 7 + 3 + (8+1) =19\]
\noindent On the other hand,
\[\left\lfloor\frac{n}{2}\sum_{v \in V(G)}\frac{c(v)-1}{c(v)}\right\rfloor = \left\lfloor\frac{8}{2}\left(7\times \frac{2}{3} + \frac{1}{2}\right)\right\rfloor = 20 \neq |E(G)|\]
\end{subcase}
\begin{subcase}
    $c(v) = c(u)$ for all $v \in V(G)$
    \newline Suppose $c(u) \geq 4$, then from \cref{2ineq2} we get;
    \begin{align*}
        1&> \frac{5}{2}\sum_{v \in V(G')}\frac{c(v)-1}{c(v)} - \frac{5}{2}\sum_{v \in V(G')}\frac{c_{G'}(v)-1}{c_{G'}(v)}\\
        &\geq\frac{5}{2}\left(5 \times \frac{3}{4}\right) - \frac{5}{2}\left(4\times\frac{2}{3}+\frac{1}{2}\right) = \frac{35}{24}
    \end{align*}
    We get a contradiction. Therefore $c(u)$ must be $3$. Thus $|N_v| = 2$ for all $v \in V(G')$. We get;
    \[|E(G)| = |E(G')| + |E(C)| + |E(G':C)| = 7 + 3 + (5\times 2) =20\]
    \noindent On the other hand,
    \[\left\lfloor\frac{n}{2}\sum_{v \in V(G)}\frac{c(v)-1}{c(v)}\right\rfloor = \left\lfloor\frac{8}{2}\left(8\times \frac{2}{3}\right)\right\rfloor = 21\neq |E(G)|\]
\end{subcase}
\noindent If $G' \cong X$, $G$ can not be an extremal graph.
\end{case}
\begin{case}
    $G' \cong Y$
    \newline Since $Y$ contains a $K_2$, clearly $c(u) \geq 2$. Note that for $2$ vertices of of $G \cong X$, $c_{G'}(v) = 2$ and for the remaining one vertex $c_{G'}(v) =1$.
    \begin{subcase}
There exists $v_0 \in V(G')$ such that $c(v_0) = c(u) -1$ and $c(v) = c(u)$ for all $v \in V(G) \setminus\{v_0\}$.
\newline Suppose $c(u) \geq 3$, then from \cref{2ineq2} we get;
\begin{align*}
    1&> \frac{3}{2}\sum_{v \in V(G')}\frac{c(v)-1}{c(v)} - \frac{3}{2}\sum_{v \in V(G')}\frac{c_{G'}(v)-1}{c_{G'}(v)}\\
    &\geq \frac{3}{2}\left(2 \times \frac{2}{3} + \frac{1}{2}\right)- \frac{3}{2}\left(2\times\frac{1}{2}+0\right) = \frac{5}{4}
\end{align*}
We get a contradiction. Therefore, $c(u)$ must be $2$. But then $c(v_0) = c(u) -1 = 1$, that is $v_0$ is isolated in $G$. Thus, we infer that $G$ a non-trivial disconnected graph other than $Y \implies G$ is not extremal (from \cref{conn}).
    \end{subcase}
    \begin{subcase}
        $c(v) = c(u)$ for all $v \in V(G)$.
        \newline Suppose $c(u) \geq 3$, then from \cref{2ineq2} we get;
        \begin{align*}
    1&> \frac{3}{2}\sum_{v \in V(G')}\frac{c(v)-1}{c(v)} - \frac{3}{2}\sum_{v \in V(G')}\frac{c_{G'}(v)-1}{c_{G'}(v)}\\
    &\geq \frac{3}{2}\left(3 \times \frac{2}{3} \right)- \frac{3}{2}\left(2\times\frac{1}{2}+0\right) = \frac{3}{2}
\end{align*}
We get a contradiction. Therefore, $c(u)$ must be $2$. Thus $|N_v| =1$ for all $v \in V(G')$. We get;
\[|E(G)| = |E(G')| + |E(C)| + |E(G':C)| = 1 + 1 + 3 =5\]
\noindent On the other hand,
\[\left\lfloor\frac{n}{2}\sum_{v \in V(G)}\frac{c(v)-1}{c(v)}\right\rfloor = \left\lfloor\frac{5}{2}\left(5\times \frac{1}{2}\right)\right\rfloor = 6 \neq |E(G)|\]
    \end{subcase}
    \noindent If $G' \cong Y$, $G$ is not an extremal graph.
\end{case}

\subsubsubsection{The if part}\label{if}
For the other direction, assume $G \in \mathcal{S} \cup \{X,Y\}$. We will show equality holds in \Cref{th:main}. \begin{itemize}
    \item If $G = X$, $|E(X)| = 7$.
\[\left\lfloor\frac{n}{2}\sum_{v \in V(X)}\frac{c(v)-1}{c(v)}\right\rfloor = \left\lfloor\frac{5}{2}\left(4\times \frac{2}{3}+\frac{1}{2}\right)\right\rfloor = 7 = |E(X)|\]
\item If $G = Y$, $|E(Y)| =1$.
\[\left\lfloor\frac{n}{2}\sum_{v \in V(Y)}\frac{c(v)-1}{c(v)}\right\rfloor = \left\lfloor\frac{3}{2}\left(2\times \frac{1}{2}+0\right)\right\rfloor = 1 = |E(Y)|\]
\item If $G \in \mathcal{S}$, let $G = T(n,r)$ for some $r$. We have $c(v) =r$ for all $v \in V(G)$. Thus from the definition of $\mathcal{S}$, in \cref{defnS}, we get;
\[|E(G)| = \left\lfloor\frac{n^2(r-1)}{2r}\right\rfloor = \left\lfloor\frac{n}{2}\left(n\times \frac{c(u)-1}{c(u)}\right)\right\rfloor = \left\lfloor\frac{n}{2}\sum_{v \in V(G)}\frac{c(v)-1}{c(v)}\right\rfloor\]
\end{itemize}

\noindent Therefore, from \ref{onlyif} and \ref{if} we get that $G$ is extremal for \Cref{th:main} if and only if $G\in \mathcal{S} \cup\{X,Y\}$.

\subsection{Localization of \Cref{wturan}}
\begin{thm}\label{wturanlocal}
  For a simple graph $G$ with $n$ vertices
\begin{equation*}
    |E(G)| \leq \frac{n}{2}\sum_{v\in V(G)}\frac{c(v)-1}{c(v)}
\end{equation*}
Equality holds if and only if $G$ is a regular Tur\'{a}n graph; that is, all the classes are of equal size. 
\end{thm}
\noindent We will closely follow the proof of \Cref{th:main} with appropriate modifications (like dropping the floor function), to develop this proof. 
\begin{proof}
 We will follow a similar induction approach, on the number of vertices, as in the proof of \Cref{th:main}.
 \newline \textit{Base Case:} The claim is trivially true when $n =1$, since both the sides are zero.
 
 \vspace{1mm}
 \noindent \textit{Induction Hypothesis: }Suppose the claim is true for all graphs with less than $|V(G)|$ vertices.

 \vspace{1mm}
 \noindent \textit{Induction Step:} Note that, the inequality in \Cref{wturanlocal} follows directly from the bound in \Cref{th:main}. It can also be derived independently by following the same steps as in the proof of \Cref{Proof of Ineq}, omitting the floor function. Now we will characterize the extremal graph class.
    \newline It is easy to verify that if $G$ is a regular Tur\'{a}n graph, that is, $G \cong T(n,r)$ where $r\mid n$, then $G$ is extremal for \Cref{wturanlocal}.
    \newline Now, suppose $G$ is extremal. Then, paralleling the arguments used in the proof of \Cref{th:main}, the following conditions must hold:
    \begin{itemize}
        \item For equality in the first inequality of \cref{ineq12}, $G'$ must be regular Tur\'{a}n graph by the induction hypothesis, where $G' = G\setminus V(C)$ and $C$ is the largest clique of $G$. Recalling our case analysis for proof of \Cref{th:main}, the cases $G' \cong X$ and $G' \cong Y$ are not relevant here. The only case we have to handle is $G' \in \mathcal{S}$. Note that if $G'$ is a regular Tur\'{a}n graph, then $G' \in \mathcal{S}$.
        \item For equality in \cref{ineq4}, without the floor function, we need \cref{eqn30} to hold, thus we must have $c(v) = c(u)$, for all $v \in V(G)$. This means \cref{sbc1.1} is not relevant.
        \begin{equation}
             c(u)\sum_{v\in V(G')}\frac{c(v) -1}{c(v)}  =  \sum_{v\in V(G')}c(v)\frac{c(v) -1}{c(v)}\label{eqn30}
        \end{equation}
        
        \item For equality in the second inequality of \cref{ineq12}, without the floor function, we need \cref{eqn31} to be true, thus we must have $c_{G'}(v) = c(v)$, for all $v \in V(G')$.  Thus only \cref{1.2.3} is possible. 
        \begin{equation}\label{eqn31}
            \frac{n-c(u)}{2}\sum_{v \in V(G')}\frac{c_{G'}(v) -1}{c_{G'}(v)}= \frac{n-c(u)}{2}\sum_{v \in V(G')}\frac{c(v) -1}{c(v)}
        \end{equation}
        \item For equality in \cref{ineq3}, we must have $|N_v| = c(v) -1$.
        \item With the previous equalities and the fact that the floor functions are omitted, \cref{ineq5} and \cref{ineq6} become equalities. Thus they are not relevant.
    \end{itemize}
    
    
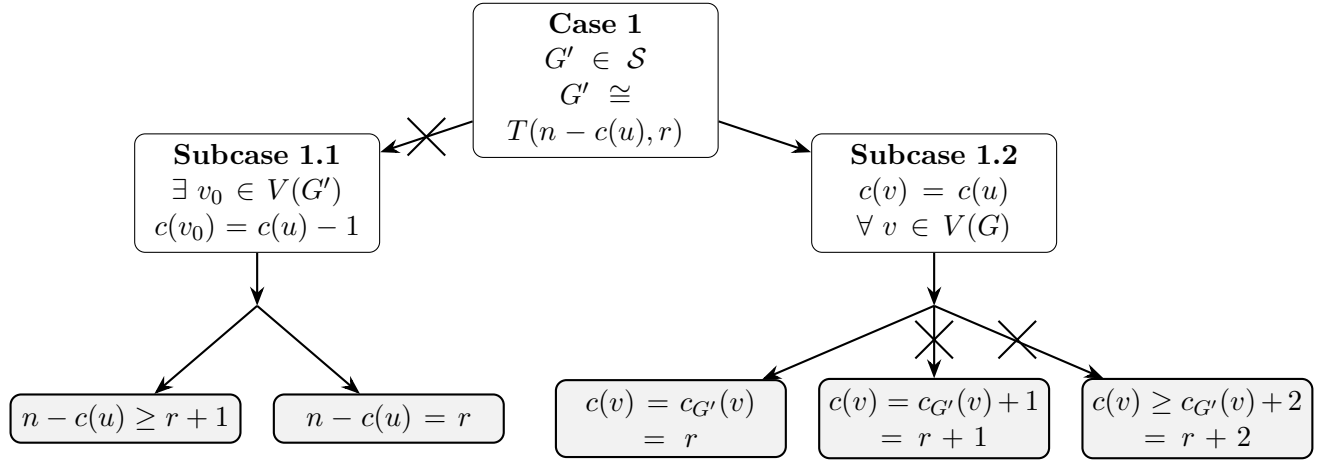
\begin{center}
\begin{tikzpicture}[
  level distance=1.5cm,
  edge from parent/.style={line},
  every node/.style={rectangle, draw, align=center, rounded corners, text width=3cm}, 
  smallbox/.style={rectangle, draw, fill=gray!10, rounded corners, align=center, text width=2.8cm},
  grow'=down,
  level 1/.style={sibling distance=90mm},
  remember picture
]

\node (root) {\textbf{Case 1}\\$G' \in \mathcal{S}$\\$G' \cong T(n-c(u) ,r)$}
  child {node (sub2) {\textbf{Subcase 1.2}\\$c(v) = c(u)$\\$\forall\ v \in V(G)$}
    child[grow'=down, level distance=1.5cm, sibling distance=35mm] {
      node[coordinate] {}
      child {node (rplus2) [smallbox] {$c(v) \geq c_{G'}(v) + 2$ \\ $= r+2$}}
      child {node (rplus1) [smallbox] {$c(v) = c_{G'}(v) + 1$ \\ $=r+1$}}
      child {node (r) [smallbox] {$c(v) = c_{G'}(v)$ \\ $=r$}}
    }
  }
  child {node (infeas) {\textbf{Subcase 1.1}\\$\exists$ $v_0 \in V(G')$\\$c(v_0) = c(u) - 1$}
    child[grow'=down, level distance=1.5cm, sibling distance=35mm] {
      node[coordinate] {}
      child {node[smallbox] {$n - c(u) = r$}}
      child {node[smallbox] {$n - c(u) \geq r+1$}}
    }
  };

\path let \p1 = (root), \p2 = (infeas) in
  coordinate (mid1) at ($(\p1)!.5!(\p2)$);
\draw[black, thick]
  ($(mid1)+(-0.25,0.25)$) -- ($(mid1)+(0.25,-0.25)$)
  ($(mid1)+(-0.25,-0.25)$) -- ($(mid1)+(0.25,0.25)$);

\path let \p1 = (sub2), \p2 = (rplus2) in
  coordinate (mid2) at ($(\p1)!.65!(\p2)+ (-1.17,0)$);
\draw[black, thick]
  ($(mid2)+(-0.25,0.25)$) -- ($(mid2)+(0.25,-0.25)$)
  ($(mid2)+(-0.25,-0.25)$) -- ($(mid2)+(0.25,0.25)$);

\path let \p1 = (sub2), \p2 = (rplus1) in
  coordinate (mid3) at ($(\p1)!.65!(\p2)$);
\draw[black, thick]
  ($(mid3)+(-0.25,0.25)$) -- ($(mid3)+(0.25,-0.25)$)
  ($(mid3)+(-0.25,-0.25)$) -- ($(mid3)+(0.25,0.25)$);

\end{tikzpicture}
\end{center}
Let the classes of $G'$ be $P_1,P_2,\dots,P_{c(u)}$. Note that all these classes are equal-sized, since $G'$ is regular Tur\'{a}n. Let $V(C)= \{u_1,u_2,\dots,u_{c(u)}\}$. Following similar arguments as in \cref{1.2.3}, let $u_i$ be the exclusive non-neighbor of the vertices of $P_i$. Thus we get $G$ to be a regular Tur\'{a}n graph with classes $P'_1,P'_2,\dots,P'_{c(u)}$, where $P'_i = P_i \cup \{u_i\}$, for all $i \in [c(u)]$.
\end{proof}
\bibliographystyle{plain}
\bibliography{references}
\end{document}